\documentclass{article}
\usepackage[T1]{fontenc}
\usepackage[mathlines]{lineno} 
\usepackage{amsmath, amsthm, amsfonts, amssymb}
\usepackage{graphicx}
\usepackage{color}
\usepackage{lineno}
\usepackage{thm-restate}
\usepackage{enumitem}
\usepackage[dvipsnames]{xcolor}
\usepackage[utf8]{inputenc}
\usepackage[english]{babel}
\usepackage{cite}
\usepackage[title]{appendix}
\usepackage{hyperref}

\usepackage{verbatim}
\usepackage{lineno} 
\usepackage{hyperref}
\usepackage{graphicx}
\usepackage{color}  
\newtheorem{thm}{Theorem}
\newtheorem{cor}[thm]{Corollary}
\newtheorem*{prop}{Proposition}
\newtheorem{prob}{Open problem}
\usepackage{bm}
 
\title{Hamiltonicity of the Complete Double Vertex Graph of some Join Graphs}
\author{Luis Manuel Rivera\thanks{Unidad Acad\'emica de Matem\'aticas, Universidad Aut\'onoma de Zacatecas, Zacatecas, Mexico. \texttt{luismanuel.rivera@gmail.com}} \and Ana Laura Trujillo-Negrete\thanks{Departamento de Matem\'aticas, Cinvestav, CDMX, Mexico. Partially supported by CONACYT (Mexico), grant 253261. \texttt{ltrujillo@math.cinvestav.mx}}
}
\date{\today}

\begin{document}
	\maketitle

\begin{abstract}

The complete double vertex graph $M_2(G)$ of $G$ is defined as the graph whose vertices are the $2$-multisubsets of $V(G)$, and two of such vertices are adjacent in $M_2(G)$ if their symmetric difference (as multisets) is a pair of adjacent vertices in $G$. In this paper we exhibit an infinite family of graphs $G$ (containing Hamiltonian and non-Hamiltonian graphs) for which $M_2(G)$ are Hamiltonian. 

\end{abstract}

\noindent
{\it Keywords:} Complete Double Vertex Graphs, Hamiltonicity, Token graphs. \\
{\it 2000 AMS Subject Classification:} 05C45, 05C76. \\


\section{Introduction}
\begin{figure}
	\centering
	\includegraphics[width=0.8\textwidth]{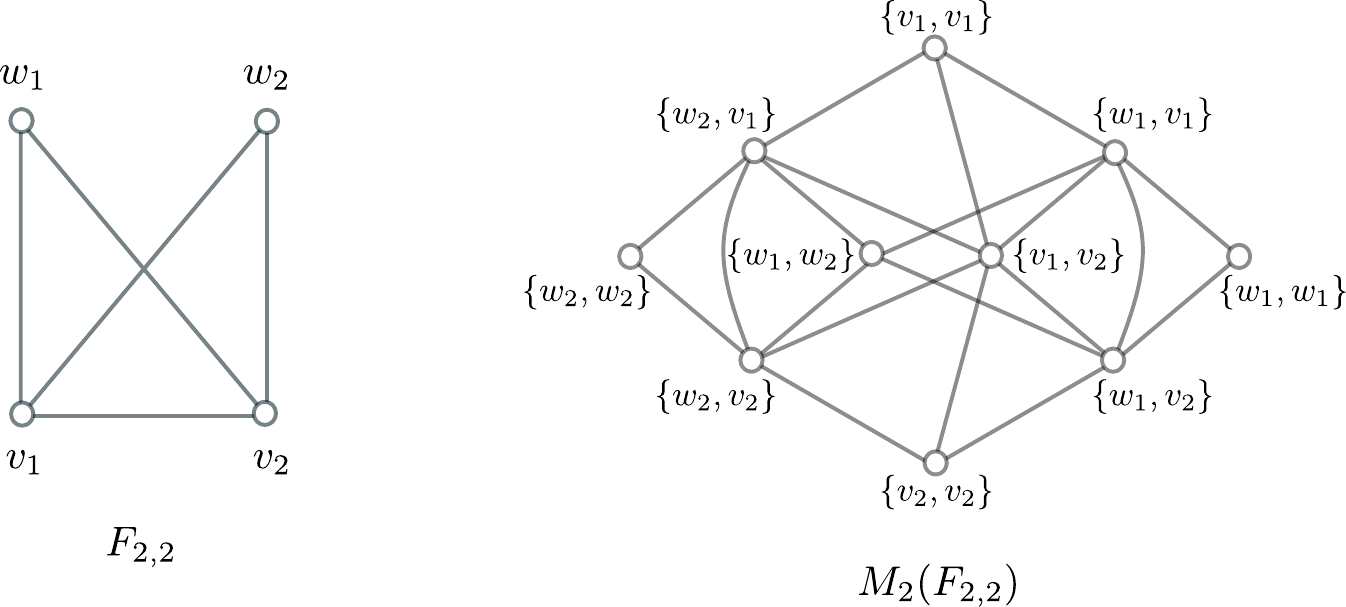}
	\caption{Graph $F_{2, 2}$ and its complete double vertex graph $M_2(F_{2, 2})$.}
	\label{fig:cdvg_example}
\end{figure}

Throughout this paper, $G$ is a simple graph of order $n \geq 2$. The {\it complete double vertex graph}  $M_2(G)$ of $G$ is the graph whose vertices are all the $2$-multisubsets of $V(G)$ where two vertices are adjacent whenever their symmetric difference (as multisubsets) is a pair of adjacent vertices in $G$. See an example in Figure~\ref{fig:cdvg_example}. This class of graphs can be seen as a generalization of the $2$-token graph $F_2(G)$ of $G$, where the vertex set consists of all the $2$-subsets of $V(G)$ and the adjacencies are defined in a similar way. It is easy to show that $G$ and $F_2(G)$ are isomorphic to a subgraph of $M_2(G)$. 

The $2$-token graph is, indeed, a particular case of the $k$-token graphs. For a graph $G$ of order $n$ and an integer $k$ with $1\leq k\leq n-1$, the \textit{$k$-token graph $F_k(G)$ of $G$} is the graph whose vertices are the $k$-subsets of $V(G)$, where two of such $k$-subsets are adjacent if their symmetric difference is a pair of adjacent vertices in $G$.

The properties of token graphs have been studied  since 1991 by various authors and with different names, see, e.g., \cite{hamiltonicity_joingraphs,alavi2,alavi3,alavi4, aude, FFHH} and, in recent years,  the study of its combinatorial properties and its connection with problems in other areas such as Coding Theory and Physics has increased, see, e.g., \cite{dealba2, token2, deepa2,  ruyanalea, fisch2, soto, leatrujillo, ouy}. However, there are only few papers about complete double vertex graphs, which were implicitly introduced by Chartrand et al.~\cite{char}. As far as we know, the only parameters of complete double vertex graphs that have been studied are: connectivity, regularity, planarity, chromatic number, Hamiltonicity and Eulerianicity. These parameters were studied by  Jacob, Goddard and Laskar~\cite{jacob} in 2007.

In this note we are interested in the Hamiltonicity of complete double vertex graphs.   Jacob, Goddard and Laskar~\cite{jacob} showed that the Hamiltonicity of $G$ does  not necessarily imply the Hamiltonicity of $M_2(G)$, for example if $n=4$ or $n \geq 6$, then $M_2(C_n)$ is not Hamiltonian. However, the authors also show that if $G'$ is the graph obtained from $C_n$ by adding an edge  between two vertices at distance two, then $M_2(G')$ is Hamiltonian. To our knowledge, these are the only known results about the Hamiltonicity problem in this class of graphs. 

In order to formulate our main result we need the following definitions. 
Given two disjoint graphs $G$ and $H$, the \textit{join graph} $G+H$ 
of $G$ and $H$ is the graph whose vertex set is $V(G)\cup V(H)$ and its edge set is 
$E(G)\cup E(H)\cup \{xy:x\in G\text{ and }y\in H\}$. The \emph{generalized fan graph}, $F_{m,n}$, or  simply the \emph{fan graph}, is the join graph $F_{m,n}=E_m +P_n$, 
where $E_m$ denotes the empty  graph on $m$ vertices ($E_m=\overline{K_m}$) and $P_n$ denotes the path  graph on $n$ vertices. In Figure 1 we show $F_{2, 2}$.

The authors of this article studied the Hamiltonicity of $k$-token graphs of fan graphs $F_{1, n}$~\cite{rive-tru}. In this note, we  continue with this line of research but for the case of the complete double vertex graph of the generalized fan graphs $F_{m, n}$, for $m\geq 1$. Our main result is the following.

\begin{thm}
	\label{thm:main1}
	The complete double vertex graph of $F_{m,n}$ is Hamiltonian if and only if $n\geq 2$ and $1\leq m \leq 2\,(n-1)$. 
\end{thm}

A direct consequence of Theorem~\ref{thm:main1} is the following result.

\begin{cor}
	\label{cor:main2}
	Let $G_1$ and $G_2$ be two graphs
	of order $m\geq 1$ and $n\geq 2$, respectively, such that $G_2$ has a Hamiltonian path. Let $G=G_1+G_2$. If  $m\leq 2(n-1)$ then $M_2(G)$ is Hamiltonian. 
\end{cor}

 Observe that Corollary~\ref{cor:main2} provides an infinite family of non-Hamiltonian graphs whose complete double vertex graphs are Hamiltonian.

\section{Proof of Theorem~\ref{thm:main1}}
\label{sec:complete}

 We introduce the following notation. 
Let $V(P_n):=\{v_1,\ldots,v_n\}$ 
and $V(E_m):=\{w_1,\ldots, w_m\}$, so   
$V(F_{m,n})=\{v_1,\ldots,v_n,w_1,\ldots,w_m\}$. 
For a path $T=a_1 a_2 \dots a_{l-1}a_{l}$, we will use $\overleftarrow{T}$ to denote the reverse path $a_la_{l-1}\dots a_2a_1$. 
As usual, for a positive integer $r$ let $[r]:=\{1,2,\ldots,r\}$.

If $n=1$ and $m\geq 1$, then $\deg(\{w_i,w_i\})=1$ for any $i\in [m]$, which implies that 
$M_2(F_{m,n})$ is not Hamiltonian, so we assume that $n\geq 2$.  

The proof is by construction. We distinguish four cases: either $m=1$, $m=2\,(n-1)$, $1<m<2\,(n-1)$ or $m>2(n-1)$. 

\begin{itemize}
	\item {\bf Case} $\boldsymbol{m=1}$

	For $1\leq i\leq n$ let 
	\[
	T_i:=\{v_i,w_1\}\{v_i,v_i\}\{v_i,v_{i+1}\}\{v_i, v_{i+2}\}\ldots \{v_i,v_n\}.
	\] 
	Note the following
	\begin{itemize}
		\item[(a)]  $T_i$ is a path in $M_2(F_{m,n})$, where $T_n=\{v_n,w_1\}\{v_n,v_n\}$;
		\item[(b)] for any $i\in [n-1]$ the final vertex of $T_i$ is adjacent to the initial vertex of $\overleftarrow{T_{i+1}}$, so the paths $T_i$ and $\overleftarrow{T_{i+1}}$ can be concatenated as $T_i\overleftarrow{T_{i+1}}$;
		\item[(c)] for any $i\in [n-1]$ the final vertex of $\overleftarrow{T_i}$ is adjacent to the initial vertex of $T_{i+1}$, so the paths $\overleftarrow{T_i}$ and $T_{i+1}$ can be concatenated as $\overleftarrow{T_i}T_{i+1}$. 
	\end{itemize}
	
	Let
	\[
	C:= \begin{cases}
		T_1\, \overleftarrow{T_2}\, T_3 \, \overleftarrow{T_4}\,\dots\, T_{n-1}\,\overleftarrow{T_n} \,\{w_1,w_1\} & \text{if $n$ is even, } \\
		\overleftarrow{T_1}\,\{w_1,w_1\}\,T_2\,\overleftarrow{T_3}\,T_4\ldots\,T_{n-1}\,\overleftarrow{T_n} & \text{if $n$ is odd.}
	\end{cases}
	\]
	We claim that $C$ is a Hamiltonian cycle in $M_2(F_{m,n})$. Suppose that $n$ is even. Statements (b) and (c) together imply that the concatenation
	$T_1\, \overleftarrow{T_2}\, T_3 \, \overleftarrow{T_4}\,\dots\, T_{n-1}\,\overleftarrow{T_n}$ produces a walk in $M_2(F_{m,n})$. Now, since $\{v_n,w_1\}$ (the final vertex of $\overleftarrow{T_n}$) is adjacent to $\{w_1,w_1\}$, and $\{w_1,w_1\}$ is adjacent to  $\{v_1,w_1\}$ (the initial vertex of $T_1$), it follows that $C$ is a closed walk in $M_2(F_{m,n})$.  Finally, note that $\big\{V(T_1), \dots, V(T_n), \{\{w_1, w_1\}\}\big\}$ is a partition of the vertex set of $M_2(F_{m,n})$. Therefore, $C$ is a Hamiltonian cycle. By similar arguments,  when $n$ is  odd we have that $C$ is a Hamiltonian cycle in $M_2(F_{m,n})$, as claimed. 
	
	Let us note here that $\{v_n,w_1\}$ and $\{v_n,v_n\}$ are adjacent in $C$ (these two vertices 
	correspond to the vertices of $\overleftarrow{T_n}$).  This observation will be useful in the following two cases.
	
	\item {\bf Case} $\boldsymbol{m=2\,(n-1).}$
	
	Let $C$ be the cycle defined in the previous case,  depending on the parity of $n$. Let 
	\[
	P_1:=\{v_n, w_1\} \xrightarrow{C} \{v_n, v_n\}
	\]
	be the path obtained from $C$ by deleting the edge between $\{v_n, w_1\}$ and $\{v_n, v_n\}$. 
	For $1< i\leq n-1$ let
	{\small
	\begin{equation}\label{eq1}
	\begin{split}
	P_i&=\{w_i,v_n\}\{w_i,w_i\}\{w_i,v_{n-1}\}\{w_i,w_1\}\{w_i,v_{n-2}\}\{w_i,w_{i+(n-2)}\}\\
	&\{w_i,v_{n-3}\}\{w_i,w_{i+(n-3)}\}\ldots \{w_i,v_2\}\{w_i,w_{i+2}\}\{w_i,v_1\}\{w_i,w_{i+1}\}. 
	\end{split}
	\end{equation}
	}
	Note that after the vertex $\{w_i, w_1\}$, the vertices of $P_i$ follow the pattern $\{w_{i}, v_j\}\{w_{i}, w_{i+j}\}$, from $j=n-2$ to $1$.
	
	For $n\leq i\leq m$ let
	{\small
	\begin{equation}\label{eq2}
	\begin{split}
	P_i&=\{w_i,v_n\}\{w_i,w_i\}\{w_i,v_{n-1}\}\{w_i,w_{i+(n-1)}\}\{w_i,v_{n-2}\}\\&
	\{w_i,w_{i+(n-2)}\}\ldots
	 \{w_i,v_2\}\{w_i,w_{i+2}\}\{w_i, v_1\}\{w_i,w_{i+1}\},
	\end{split}
	\end{equation}}
	where the sums are taken mod $m$ and with the convention that $m\pmod{m}=m$. In this case, after $\{w_i, w_i\}$, the vertices in $P_i$  follow the pattern $\{w_{i}, v_j\}\{w_{i}, w_{i+j}\}$, from $j=n-1$ to $1$.
	We claim that 
		\[
	P:=P_1\,P_2\,\ldots\,P_m\{v_{n}, w_{1}\}
	\]
	is a Hamiltonian cycle in $M_2(F_{m,n})$. 
	
	Note the following 
	\begin{enumerate}
		\item the final vertex of $P_1$ is $\{v_n,v_n\}$ while the initial vertex of $P_2$ is $\{w_2,v_n\}$, and 
		these two vertices are adjacent in $M_2(F_{m,n})$; 
		\item for $1<i<2(n-1)$, the final vertex
		of $P_i$ is $\{w_i,w_{i+1}\}$ while the initial vertex of $P_{i+1}$ is $\{w_{i+1},v_n\}$, and also
		these two vertices are adjacent in $M_2(F_{m,n})$;
		\item the final vertex of $P_m$ is $\{w_m,w_1\}$ while
		the initial vertex of $P_1$ is $\{v_n,w_1\}$, and these two vertices are adjacent in  $M_2(F_{m,n})$.
	\end{enumerate}   
	 These three observations together imply that $P$ is a closed walk in $M_2(F_{m,n})$. 
	Now we show that $\{V(P_1), \dots, V(P_m)\}$ is a partition of $V(M_2(F_{m,n}))$. 
	Notice that any vertex in $V(F_{m,n})$  satisfies exactly one of the following conditions:  
	\begin{itemize}
		\item The vertices of type $\{v_i,v_j\}$   belong to $P_1$	for any $i,j\in [n]$.
	\item The vertices of type $\{w_i,v_j\}$   belong to $P_i$ for any $i\in [m]$ and $j\in [n]$. 
	\item The vertices of type   $\{w_i,w_i\}$ belong to $P_i$ for any $i\in [m]$. 
	\item Let $H:=\{\{w_i,w_j\}\in V(M_2(F_{m,n})) \colon i,j\in [m], i\neq j\}$. Note that \[|H|=\binom{m}{2}=\frac{m(m-1)}{2}=(n-1)(m-1),\] and furthermore, $|V(P_\ell)\cap H|=n-1$ for each $\ell\in \{2,\dots,m\}$.  Moreover, for a vertex $\{w_i,w_j\}$, with $i,j\in [m]$ and $i<j$, we claim that
	\[\{w_i,w_j\}\in \begin{cases}
		P_j & \text{if $i=1$; }\\
		P_i & \text{if $i>1$ and $|j-i|\le n-2$; } \\
		P_j & \text{if $i>1$ and $|j-i|>n-2$.}
	\end{cases}\]
First suppose $i=1$. If $1<j\le n-1$, then Equation (\ref{eq1}) implies that $\{w_1,w_j\}\in P_j$, and if $n\le j\le m$, then $w_1\in \{w_{j+1},w_{j+2},\dots,w_{j+(n-1)}\}$ and hence Equation (\ref{eq2}) implies that $\{w_1,w_j\}\in P_j$. 
Suppose now that $i>1$. If $j-i\le n-2$, then $w_j\in \{w_{i+1},w_{i+2},\dots,w_{i+(n-2)}\}$, and so by Equation (\ref{eq1}) it follows that $\{w_i,w_j\}\in P_i$. If $j-i>n-2$ then $w_i\in \{w_{j+1},w_{j+2},\dots,w_{j+(n-1)}\}$, and so by Equation (\ref{eq2}) we have that $\{w_i,w_j\}\in P_j$. 
	\end{itemize}
	
	Thus,  $P$ is a Hamiltonian cycle in $M_2(F_{m,n})$.
	
	\item {\bf Case}  $\boldsymbol{1<m<2\,(n-1).}$
	
	We consider again the paths $P_1,\ldots,P_m$ defined in the previous case 
	with a slight modification: 
	\begin{itemize}
		\item $P_1'=P_1$;
		\item for $i \in \{2, \dots, m-1\}$, let $P'_i$ be the path obtained from $P_i$ by deleting the vertices of type $\{w_i, w_j\}$, for each $j>m$;
		\item let $P'_m$ be the path obtained from $P_m$ by first interchanging $\{w_m, w_{m+1}\}$ and $\{w_m, w_1\}$ from their current positions in $P_m$, and then deleting the vertices of type $\{w_m, w_j\}$, for every $j>m$.
	\end{itemize}
	We have  that $P'_1,\ldots, P'_m$ are, indeed, disjoint paths in $M_2(F_{m,n})$, and that
	$P'_i$ has the same initial and final vertices  as $P_i$ for $i \neq m$, so the concatenation 
	\[
	P':=P'_1\,\ldots\,P'_m\{v_n, w_1\}
	\]
	 corresponds to a cycle in $M_2(F_{m,n})$. It is an easy exercise 
	to show that $P'$ is, in fact, a Hamiltonian cycle of $M_2(F_{m,n})$. 
	
	\item {\bf Case}  $\boldsymbol{m>2\,(n-1).}$

In this case we show that $M_2(F_{m,n})$ is not Hamiltonian by using the following result  stated in  West's book~\cite{west}. 
    \begin{prop}[Prop. 7.2.3, \cite{west}]
    	If $G$ has a Hamiltonian cycle, then for each nonempty set $S\subseteq V(G)$, the graph $G-S$ has at most $|S|$ connected components. 
    \end{prop}	
   So, for our purpose, we are going 
	to exhibit a subset $S\subset V(M_2(F_{m,n}))$
	such that the graph $M_2(F_{m,n})-S$ has at most $|S|$ connected components.
	Let 
\[
\begin{aligned}
	S:=&\{\{w_i,v_j\}\in V(M_2(F_{m,n})) : i\in [m]\text{ and }j\in [n]\}, \\
	T:=&\{\{w_i,w_j\}\in V(M_2(F_{m,n})) : i,j\in [m]\}, \\
	R:=&\{\{v_i,v_j\}\in V(M_2(F_{m,n})) : i,j\in [n]\}.
\end{aligned} 
\]
The set $\{S, T, R\}$  is a partition of $V(M_2(F_{m,n}))$. Since any vertex in $T$ has its neighbors in $S$, 
	the subgraph induced by $T$ in $M_2(F_{m,n})-S$ is isomorphic to $\overline{K_{\binom{m+1}{2}}}$. 
	On the other hand, the subgraph induced by $R$ is a component of $M_2(F_{m,n})-S$.  
	Since $|S|=m\,n$, $|T|=\binom{m+1}{2}$ and $m>2\,(n-1)$,  we have that the number of connected components of $M_2(F_{m,n})-S$ is
	\[
 |T|+1=\binom{m+1}{2}+1> m\,n=|S|.\]
	This completes the proof of Theorem~\ref{thm:main1}. 
\end{itemize}

\section{Open problems}

There is a natural generalization of the complete double vertex graphs: Let $G$ be a simple graph of $n$ vertices and let $k$ be an integer with $1\leq k\leq n-1$. The \textit{$k$-multiset graph $M_k(G)$ of $G$} is the graph whose vertices are the $k$-multisubsets of $V(G)$, with two vertices being adjacent in $M_k(G)$ if their symmetric difference (as multisets) is a pair of adjacent vertices in $G$. This generalization gives rise to several open problems. 

\begin{prob}
	To  establish whether or not the Hamiltonicity of $G$ implies the Hamiltonicity of $M_k(G)$  for $k>2$. 
\end{prob}

In this paper we provided an infinite family of non-Hamiltonian graphs with Hamiltonian complete double vertex graphs. A natural open problem is to study if this result can be generalized.  

\begin{prob}
	Given  an integer $k>2$, to establish whether or not there exist non-Hamiltonian graphs with Hamiltonian $k$-multiset graphs.  
\end{prob}


\newpage 

\end{document}